\newcommand{\rest}{\mathord{\restriction}}
\renewcommand{\phi}{\varphi}

\newcommand{\sat}{\models}
\newcommand{\su}{\subseteq}
\renewcommand{\a}{\alpha}
\renewcommand{\b}{\beta}
\newcommand{\g}{\gamma}
\newcommand{\e}{\varepsilon}

\renewcommand{\d}{\delta}
\renewcommand{\l}{\lambda}
\renewcommand{\k}{\kappa}
\newcommand{\z}{\zeta}

\newcommand{\om}{\omega}

\newcommand{\lng}{\langle}
\newcommand{\rng}{\rangle}
\newcommand{\ov}{\overline}
\newcommand{\sm}{\setminus}

\newcommand{\dom}{{\operatorname {dom}}}
\newcommand{\cf}{{\operatorname {cf}}}

\newcommand{\otp}{{\operatorname {otp}}}
\newcommand{\ran}{{\operatorname  {ran}}}

\newcommand{\id}{{\text {id}}}

\newcommand{\imply}{\Rightarrow}

\newcommand{\conc}{\widehat{\thickspace}}

\newcommand{\N}{\mathbb N}
\newcommand{\R}{\mathbb R}
\newcommand{\Q}{\mathbb Q}

\documentclass[11pt,reqno]{amsart}

 \usepackage{amstext} \usepackage{amssymb}

\newtheorem{lemma}{Lemma}
\newtheorem{theorem}[lemma]{Theorem}

 \newtheorem{definition}[lemma]{Definition}
 \newtheorem{corollary}[lemma]{Corollary}
 
 \newtheorem{problem}[lemma]{Problem} \newtheorem{claim}[lemma]{Claim}
 \newtheorem{fact}[lemma]{Fact}

\usepackage{graphics}
\renewcommand{\int}{\operatorname{\rm int}}

\newcommand{\U}{\mathbb U}
\newcommand{\A}{\mathbb A}

\title{Almost isometric Embeddings of metric spaces}

\author{Menachem Kojman${}^1$} \address{Department of Mathematics,
Ben-Gurion University of the Negev, Beer-Sheva, Israel}

\email{kojman@math.bgu.ac.il}

\thanks{${}^1$ Research partially
supported by an Israeli Science foundation grant no. 177/01}

\author{Saharon Shelah${}^2$}
\address{Institute of Mathematics, Hebrew University of
        Jerusalem, Israel \\ \textit{and} \\ Department of Mathematics Rutgers
        University, New-Brunswick} 

\email{shelah@math.huji.ac.il}

\thanks{Research supported by the United States-Israel Binational Science
   Foundation. Publication 827.}

\begin{document}
\maketitle

\begin{abstract}
We investigate a relations of \emph{almost isometric embedding} and
\emph{almost isometry} between metric spaces and prove that with
respect to these relations:
\begin{enumerate}
\item There is a countable universal metric space.
\item There may exist fewer than continuum separable metric spaces
  on $\aleph_1$ so that every separable metric space is almost
  isometrically embedded into one of them when the continuum
  hypothesis fails.
\item There is no collection of fewer than continuum metric spaces of
  cardinality $\aleph_2$ so that every ultra-metric space of
  cardinality $\aleph_2$ is almost isometrically embedded into one of
  them if $\aleph_2<2^{\aleph_0}$.
\end{enumerate}

We also prove that various spaces $X$ satisfy that if a space $X$ is
almost isometric to $X$ than $Y$ is isometric to $X$.

\end{abstract}

\section{Introduction}
We this paper we investigate a relation  between metric spaces
that we call  ``almost isometric embedding'', and the notion of
similarity associated with it,  ``almost isometry''. The notion
of almost isometric embedding is a weakening of the notion of
isometric embedding, with respect to which there exists a
countable universal metric space which is unique up to almost
isometry. On the other hand, almost isometry is a sufficiently
strong notion to allow many important metric spaces to maintain
their isometric identity: any space which is almost isometric to,
e.g., $\R^d$ is in fact isometric to $\R^d$.

We begin by examining two properties of separable metric spaces:
almost isometric uniqueness of a countable dense set and almost
isometric uniqueness of the whole space. It turns out that both
properties are satisfied by many well-known metric spaces. This
is done in Section \ref{unique}.

 In the rest of the paper we  investigate whether analogs of
well-known set-theoretic and model-theoretic results about
embeddability in various categories (see
\cite{baumgartner,Sh:100, mekler, KjSh:409, rep})
 hold for almost isometric embeddability in the classes of separable
  and of not necessarily separable metric spaces. Forcing is used
  only in Section \ref{forcing} and uncountable combinatorics
  appears in Section \ref{club}. Other than that, all proofs in
  the paper are elementary.

Baumgartner has shown long ago the consistency of ``all
$\aleph_1$-dense sets of $\R$ are order isomorphic'', all natural
analogs of this statement are false; in Section \ref{pairwise} it
is seen that the  statement ``all $\aleph_1$-dense subsets of $X$
are almost isometric,'' is, in contrast, false in a strong sense
for $X=\R$, $X=\R^d$ and $X=\U$, Uryson's universal metric space.
In each of these spaces there are $2^{\aleph_0}$ pairwise
incomparable subspaces with respect to almost isometric
embedding.

On the other hand, Section \ref{forcing} shows that it is
consistent that fewer than continuum $\aleph_1$-dense subsets of
$\U$ suffice to almost isometrically embed every such set, which
is a partial analog to what was proved in \cite{Sh:100,
Sh:gr1,Sh:gr2} for linear orders and graphs.

Finally, Section \ref{club} handles the category of not
necessarily separable spaces. We prove that the relation of almost
isometric embedding between metric spaces of regular cardinality
$\aleph_2$ or higher admits a representation as set inclusion
over sets of reals, similarly to what is known for linear
orderings, models of stable-unsuperstable theories, certain
groups and certain infinite graphs \cite{KjSh:409, stuss, groups,
rep, DzSh}. A consequence of this is that, in contrast to the
separable case, it is impossible to have fewer than continuum
metric spaces of cardinality $\aleph_2$ so that every metric
space on $\aleph_2$ is almost isometrically embedded into one of
them if the continuum is larger than $\aleph_2$.

The results in Sections \ref{unique} and \ref{club} were proved
by the first author and the results in Sections \ref{pairwise} and
\ref{forcing} were proved by the second author.

\section{Basic definitions and some preliminaries}
Cantor proved that the order type of $(\Q,<)$ is characterized
among all countable ordertypes as being dense and with no
end-points. Thus, any countable dense subset of $\R$ is order
isomorphic to $\Q$. Consider countable dense subsets of $\R$ with
the usual  metric. Not all of them are isometric, as $\Q$ and
$\pi \Q$ show.  In fact there are $2^{\aleph_0}$ countable dense
subsets of $\R$ no two of which are comparable with respect to
isometric embedding, since the set of distances which occur in a
metric space is an isometry invariant which is preserved under
isometric embeddings.

We introduce now relations of similarity and embeddability which are
quite close to isometry and isometric embedding, but with respect to
which the set of distances in a space is not preserved.

\begin{definition}
\begin{enumerate}
\item A map $f:X\to Y$, between  metric spaces satisfies the Lipschitz
  condition with constant $\l>0$ if for all $x_1,x_2\in X$ it holds that
  $d(f(x_1), f(x_2))< \l d(x_1,x_2)$.
\item Two metric spaces $X$ and $Y$ are \emph{almost isometric} if for each
$\l>1$ there is a homeomorphism $f:X\to Y$ such that $f$ and $f^{-1}$
satisfy the Lipschitz condition with constant $\l$.
\item $X$ is \emph{almost isometrically embedded} in $Y$ if for all
  $\l>1$ there is an injection  $f:X\to Y$ so that $f$ and $f^{-1}$
  satisfy the Lipschitz condition with constant $\l$.
\end{enumerate}
\end{definition}

Let us call, for simplicity, an injection $f$ so that $f$ and $f^{-1}$
satisfy the $\l$-Lipschitz condition, $\l$-bi-Lipchitz. Observe that
we use a strict inequality in the definition of the Lipschitz
condition.  We also note that $X$ and $Y$ are almost isometric if and
only if the \emph{Lipschitz distance} between $X$ and $Y$
is $0$. The Lipschitz distance is a well known semi-metric on metric
spaces (see \cite{Gromov}).

It is important to notice that one does not require in (2) that the
injections $f$ for different $\l$ have the same range.

The graph of the function $f(x)=1/x$ for $x>1$ and the ray
$\{(x,0):x>0\}$ are each almost isometrically embedded in the
other but are not almost isometric. In the infinite dimensional
Hilbert space it is not hard to find two closed subsets which are
almost isometric but not isometric as follows: fix an orthonormal
basis $\{v_n:n\in \N\}$ and fix a partition $\Q=A_1\cup A_2$ to
two (disjoint) dense sets and let $A_i=\{ r^i_n:n\in \N\}$ for
$i=1,2$. Let $X_i=\bigcup_{n\in \N} [0,q_n v_n]$. Now $X_1,X_2$
are closed (and connected) subsets of the Hilbert space which are
not isometric, but are almost isometric because  $A_1,A_2$ are.

Almost isometry and almost isometric embedding can be viewed as
isomorphisms and monomorphisms of a category as follows. Let $\mathcal
M$ be the category in which the objects are metric spaces and the
morphisms are defined as follows: a morphism from $A$ to $B$ is a
sequence $\vec f=\lng f_n:n\in n\rng$ where for each $n$, $f_n:A\to B$
satisfies the Lipchitz condition with a constant $\l_n$ so that
$\lim_n\l_n=1$. The identity morphism $\vec \id_A$ is the constant
sequence $\lng \id_A:n\in \N\rng$ and the composition law is $\vec
g\circ \vec f=\lng g_n\circ f_n:n\in \N\rng$.

A morphism $\ov f$ in this category is invertible if and only if each
$f_n$ is invertible and its inverse satisfies the Lipchitz condition
for some $\theta_n$ so that $\lim_n\theta_n=1$, or, equivalently,
$\vec f\in \hom(A,B)$ is an isomorphism if and only if $f_n$ satisfies
a bi-Lipschitz condition with a constant $\l_n$ so that $\lim
_n\l_n=1$. Thus, $A$ and $B$ are isomorphic in this category if and
only if for all $\l>1$ there is a $\l$-bi-Lipchitz homeomorphism
between $A$ and $B$.

\medskip
In Section \ref{unique} and in Section \ref{forcing} below we shall
use the following two simple facts:

\begin{fact}\label{amalgamation}
Suppose $X$ is a nonempty finite set, $E\su X^2$ symmetric and
reflexive, and $(X,E)$ is a connected graph. Suppose that $f:E\to
\R^+$ is a symmetric function such that $f(x,y)=0\iff x=y$ for all
$(x,y)\in E$.  Let $\sum_{i<d}f(x_i,x_{i+1})$ be the
\emph{length} of a path $(x_0,x_1,\dots,x_d)$ in $(X,E)$. Define
$d(x,y)$ as the length of the shortest path from $x$ to $y$ in
$(X,E)$. Then $d$ is a metric on $X$. If, furthermore, for all
$(x,y)\in \dom f$ it holds that $(x,y)$ is the shortest path from
$x$ to $y$, then $d$ extends $f$.
\end{fact}

\begin{fact}\label{sp} Suppose $X=\{x_0,\dots,x_{n-1},x_n\}$ and
  $Y=\{y_0,\dots,y_{n-1}\}$ are metric spaces and that the mapping
  $x_i\mapsto y_i$ for all $i<n$ is $\theta$-bi-Lipschitz,
  $\theta>1$. Then there is a metric extension of  $Y\cup
  \{y_n\}$ of $Y$,  $y_n\notin Y$, so that the extended mapping
  $x_i\mapsto y_i$ for $i\le n$ is $\theta$-bi-Lipschitz.
\end{fact}

\begin{proof}
Fix $1<\l<\theta$ so that the mapping $x_i\mapsto y_i$ is
$\l$-bi-Lipschitz. Add a new point $y_n$ to $Y$ and define
$d^*(y_n,y_i)=\l d(x_n,x_i)$ for all $i<n$, $d^*(y_n,y_n)=0$. Let $d'$
on $Y\cup \{y_n\}$ be the shortest path metric obtained from $d\cup
d^*$.

Let us verify that $d'$ extends the given metric $d$ on
$Y$. $d(y_i,y_j)\le \l d(x_i,x_j)\le \l(d(x_n,x_i) + d(x_n,x_j)) =
d^*(y_n,y_i) + d^*(y_n,y_j)$, so $(y_i,y_j)$ is the shortest path from
$y_i$ to $y_i$.

Now let us verify that extending the mapping by $x_n\mapsto y_n$
yields a $\l$-bi-Lipschitz map.  The path $(y_n,y_i)$ has length
$\l d(x_n,x_i)$ with respect to $d\cup d^*$, hence the shortest
path cannot be longer, and $d'(y_n,y_i)\le \l d(x_n,x_i)$.
Suppose now the shortest path from $y_n$ to $y_i$ is
$(y_n,y_j,y_i)$. Then $d'(y_j,y_i)\ge d(x_j,x_i)/\l$ and
certainly $d^*(y_n,y_j) > d(x_n,x_j)$, so $d'(y_n,y_i)\ge
(d(x_n,x_j)+d(x_j,x_i))/\l \ge d(x_n,x_i)/\l$ as required.
\end{proof}

\begin{definition}
Let $LAut(X)$ be the group of all auto-homeomorphisms of $X$ which
are $\l$-bi-Lipschitz for some $\l>1$. Let $LAut_\l(X)=\{f\in LAut(X):
f \text{ is $\l$-bi-Lipschitz}\}$.

$X$ is \emph{almost ultrahomogeneous} if every finite
$\l$-bi-Lipschitz map $f:A\to B$ between finite subsets of $X$
extends to a $\l$-bi-Lipschitz autohomeomorphism.
\end{definition}

\section{Almost-isometry  uniqueness and countable dense sets}
\label{unique}

\begin{definition} A metric space $X$ is \emph{almost-isometry unique}
  if every metric space $Y$ which is almost isometric to $X$ is
  isometric to $X$.
\end{definition}

In this section we shall prove that various metric spaces are
almost-isometry unique and prove that the Uryson space $\U$  has a
unique countable dense set up to almost isometry.

\begin{theorem}Suppose $X$ satisfies:
\begin{enumerate}
\item all closed bounded balls in $X$ are compact;
\item there is $x_0\in X$ and $r>0$ so that for all $y\in X$ and all
$\l>1$ there is a $\l$-bi-Lipschitz auto-homeomorphism of $X$ so that
$d(f(y),x_0)<r$.
\end{enumerate}
Then $X$ is almost-isometry unique.
\end{theorem}

\begin{proof}
 Suppose $X$, $x_0\in X$ and $r>0$ are as stated, and suppose $Y$ is a
 metric space and $f_n:Y\to X$ is a $(1+1/n)$-bi-Lipschitz
 homeomorphism for all $n>0$. Fix some $y_0\in Y$. By following each
 $f_n$ by a bi-Lipschitz autohomeomorphism of $X$, we may assume, by
 condition 1, that $d(f_n(y_0),x_0)<r$ for all $n$.

 Condition 1 implies that $X$ is separable and complete, and since $Y$
 is homeomorphic to $X$, $Y$ is also separable. Fix a countable dense
 set $A\su Y$. Since for each $a\in A$ it holds that $d(f_n(a), x_0)$
 is bounded by some $L_a$ for all $n$, condition 2 implies that there
 is a converging subsequence $\lng f_n(a): n\in D_a\rng$ and, since
 $A$ is countable, diagonalization allows us to assume that $f_n(a)$
 converges for every $a\in A$ to a point we denote by $f(a)$. The
 function $f$ we defined on $A$ is clearly an isometry, and hence can
 be extended to an isometry on $Y$. It can be verified that $f_n(y)$
 converges pointwise to $f(y)$ for all $y\in Y$. Since each $f_n$ is
 onto $X$, necessarily also $f$ is onto $X$.  Thus $X$ is isometric to
 $Y$.
\end{proof}

\begin{corollary}
For each $d\in \N$, $\R^d$ and $\mathbb H^d$ are almost-isometry unique.
\end{corollary}

\begin{theorem}[Hrusak, Zamora-Aviles \cite{HZ}] Any two countable dense subsets
  of $\R^n$ are almost isometric. Any two countable dense subsets of
  the separable infinite dimensional Hilbert space are almost
  isometric.
\end{theorem}


If one regards a fixed $\R^d$ as a universe of metric spaces,
namely considers only the subsets of $\R^d$, then among the
countable ones there is a universal element with respect to
almost isometric embedding, which is unique up to almost isometry:

\begin{corollary}Every dense subset of $\R^d$ is almost-isometry
  universal in the class of countable subspaces of $\R^d$.
\end{corollary}

\begin{proof}
Let $B\su \R^d$ be countable and extend $B$ to a countable dense
$A'\su \R^d$. Since $A'$ and $A$ are almost isometric, $B$ is
almost isometrically embedable into $A$.
\end{proof}

\subsection{The Uryson space.}

Uryson's universal separable metric space $\U$ is characterized up to
isometry by separability and the following property:

\begin{definition}[\textbf{Extension property}] A metric space $X$
  satisfies the \emph{extension property} if for every finite
$F=\{x_0,\dots,x_{n-1},x_n\}$, every isometry
$f:\{x_0,\dots,x_{n-1}\}\to X$ can be extended to an isometry $\hat
f:F\to X$.
\end{definition}

Separability together with the extension property easily imply the
isometric uniqueness of $\U$ as well as the fact that Every separable
metric space is isometric to a subspace of $\U$ and that $\U$ is
ultrahomogeneous, namely, every isometry between finite subspaces of
$\U$ extends to an auto-isometry of $\U$. This property of $\U$ was
recently used to determine the Borel complexity of the isometry
relation on polish spaces \cite{CGK,GK}.

\begin{definition}[\textbf{Almost extension property}] A metric space $X$
  satisfied the \emph{almost extension property} if  for every finite space
$F=\{x_0,\dots,x_{n-1},x_n\}$ and $\l>1$, every $\l$-bi-Lipschitz
$f:\{x_0,\dots,x_{n-1}\}\to X$ can be extended to a $\l$-bi-Lipschitz
$\hat f:F\to X$.
\end{definition}

\begin{claim}
Suppose that $A\su \U$ is dense in $\U$. Then $A$ satisfies the
almost extension property.
\end{claim}

\begin{proof}
Suppose $f:\{x_1,\dots,x_{n-1}\}\to A$ is $\l$-bi-Lipschitz. By Lemma
\ref{sp} there is a metric extension $\ran f\cup \{y_n\}$ so that
$f\cup \{(x_n,y_n)\}$ is $\l$-bi-Lipschitz. By the extension property
of  $\U$, we may assume that $y_n\in \U$. Now replace $y_n$ by a
sufficiently close $y'_n\in A$ so that $f\cup\{(x_n,y'_n)\}$ is
$\l$-bi-Lipschitz.
\end{proof}

A standard back and forth argument shows:

\begin{fact}Any two  countable metric spaces that satisfy
  the almost extension property are almost isometric.
\end{fact}

Therefore we have proved:

\begin{theorem}\label{cdu} Any two countable dense subsets of $\U$ are almost
  isometric.
\end{theorem}

A \emph{type} $p$ over a metric space $X$ is a function $p:X\to \R^+$
so that in some metric extension $X\cup \{y\}$ it holds that
$d(y,x)=p(x)$ for all $x\in X$. A point $y\in Y$ \emph{realizes} a
type $p$ over a subset $X\su Y$ if $p(x)=d(y,x)$ for all $x\in X$.
The extension property of $\U$ is equivalent to the property that
every type over a finite subset of $\U$ is realized in $\U$.

\begin{theorem}[Uryson] If a countable metric space $A$ satisfies the
  almost extension property then its completion $\bar A$ satisfies the
  extension property and is therefore isometric to $\U$.
\end{theorem}

\begin{proof} Let $X\su \bar A$ be a
  finite subset and let $p$ a metric type over $X$. Given $\e>0$, find
  $\l>1$ sufficiently close to $1$ so that for all $x\in X$ it holds
  that $\l p(x)- p(x)<\e/2$ and $p(x)-p(x)/\l<\e/2$ and find, for each
  $x\in X$ some $x'\in A$ with $d(x',x)<\e/2$ and sufficiently small
  so that the map $x\mapsto x'$ is $\l$-bi-Lipschitz. By the almost
  extension property of $A$ there is some $y\in A$ so that
  $d(y,x)<\e$. Thus we have shown that for all finite $X\su \bar A$
  and type $p$ over $X$, for every $\e>0$ there is some point $y\in A$
  so that $|d(y,x)-p(x)|<\e$ for all $x\in X$.

Suppose now that $X\su \bar A$ is finite and that $p$ is some type
  over $X$. Suppose $\e>0$ is small, and that $y\in A$ satisfies that
  $|d(y,x)-p(x)|<\e$ for all $x\in X$. Extend the type $p$ to $X\cup
  y$ by putting $p(y)=2\e$ (since $|d(y,x)-p(x)|<\e$, this is indeed a
  type). Using the previous fact, find $z\in A$ that realizes $p$ up
  to $\e/100$ and satisfies $d(y,z)<2\e$.

Iterating the previous paragraph, one gets a Cauchy sequence
$(y_n)_n\su A$ so that for all $x\in X$ it holds $d(y_n,x)\to
p(x)$. The limit of the sequence satisfies $p$ in $\bar A$.
\end{proof}

We now have:

\begin{fact}
A countable metric space is isometric to a dense subset of $\U$ if and
only if it satisfies the almost extension property.
\end{fact}

\begin{theorem}
The Uryson space $\U$ is almost isometry unique.
\end{theorem}

\begin{proof}
Suppose $X$ is almost isometric to $\U$. Fix a countable dense $A\su
\U$ and a countable dense $B\su X$. For every $\l>1$, $B$ is
$\l$-bi-Lipschitz homeomorphic to \emph{some} countable dense subset
of $\U$, so by  Theorem \ref{cdu}  it is $\l^2$-bi-Lipschitz
homeomorphic to $A$. So $A$ and $B$ are almost isometric. This shows
that $B$ has the almost extension property. By Uryson's theorem, $\bar
B=Y$ is isometric to $\U$.
\end{proof}

Let us construct now, for completeness of presentation, some countable
dense subset of $\U$.

Let $\lng A_n:n\in \N\rng$ be an increasing sequence of finite metric
spaces so that:
\begin{enumerate}
\item all distances in $A_n$ are rational numbers.
\item for every function $p:A_n\to \Q^n$ which satisfies the triangle
inequality ($p(x_1)+d(x_1,x_2)\ge p(x_2)$ and $p(x_1)+p(x_2)\ge
d(x_1,x_2)$ for all $x_1,x_2\in A_n$) \emph{and} satisfies that
$p(x)\le n+1$ is a rational with denominator $\le n+1$ there is $y\in
A_{n+1}$ so that $p(x)=d(y,x)$ for all $x\in A_n$.
\end{enumerate}

Such a sequence obviously exits. $A_0$ can be taken as a singleton. To
obtain $A_{n+1}$ from $A_n$ one adds, for each of the finitely many
distance functions $p$ as above a new point that realizes $p$, and
then sets the distance $d(y_1,y_1)$ between two new points to be
$\min\{d(y_1,x)+d(y_2,x):x\in A_n\}$.
Let $\A:=\bigcup_n A_n$. To see that $\A$ satisfies the almost
extension property, one only needs to verify that every type $p$ over
a finite rational $X$ can be arbitrarily approximated by a rational
type (we leave that to the reader).

This construction, also due to Uryson, shows that the Uryson space has
a dense rational subspace. (Another construction of $\U$ can be found
in \cite{Gromov}). This is a natural place to recall:

\begin{problem}[Erd\H os]Is there a dense rational subspace of $\R^2$?
\end{problem}

\begin{problem}Is it true that for a separable and homogeneous complete metric
  space all countable dense subsets are almost isometric if and only
  if the space is almost-isometry unique?
\end{problem}

\section{Almost-isometric embeeddability between $\aleph_1$-dense sets}
\label{pairwise}

\begin{definition} A subset $A$ of a metric space $X$ is
  \emph{$\aleph_1$-dense} if for every nonempty open ball $u\su X$ it
  holds that $|A\cup u|=\aleph_1$.
\end{definition}

Among the early achievement of the technique of forcing a place of
honor is occupied by Baumgartner's proof of the consistency of
``all $\aleph_1$-dense sets of reals are order isomorphic"
\cite{baumgartner}.  Early on, Sierpinski proved that there are
$2^{2^{\aleph_0}}$ order-incomparable continuum-dense subsets of
$\R$, hence Baumgrtner's consistency result necessitates the
failure of CH.

Today it is known that Baumgartner's result follows from forcing
axioms like PFA and Martin's Maximum, and also from Woodin's
axiom $(*)$
 \cite{woodin}.

 In our context, one can inquire about the consistency
  two natural analogs of the statement
 whose consistency was established by
Baumgartner. First, is it consistent that all $\aleph_1$-dense
subsets of $\R$ are almost isometric? Since every bi-Lipschitz
homeomorphism between two dense subset of $\R$ is either order
preserving or order inverting, this statement strengthens a
slight relaxation of Baugartner's consistent statement. Second,
since $\U$ in the category of metric spaces has the role $\R$ has
in the category of separable linear orders (it is the universal
separable object), is it consistent that all $\aleph_1$-dense
subsets of $\U$ are almost isometric?

The answer to both questions is negative.

\subsection{Perfect incomparable subsets in the cantor space}
For every infinite $A\su 2^\N$ let $T_A\su 2^{<\N}$ be defined
(inductively) as the tree $T$ which contains the empty sequence
and for every $\eta\in T$ contains $\eta\conc 0, \eta\conc 1$ if
$|\eta|\in A$, and contains only $\eta\conc 0$ if $|\eta|\notin
A$. Let $D_A$ be the set of all infinite branches through $T_A$.
The set of positive distances occurring in $D_A$, namely
$\{d(x,y):x,y\in D_A,x\not=y\}$ is equal to $\{1/2^n:n\in A\}$.

Let $d_3$ be the metric on $2^\N$ defined by
$d_3(\eta,\nu):=1/3^{\Delta)\eta,\nu)}$. Observe that the natural
isomorphism between $2^\N$ and the standard "middle-third" cantor set
is a bi-Lipschitz map when $2^\N$ is taken with $d_3$.

For infinite sets $A,B\su N$ let us define the following
condition:

\medskip

$(*)$ \textit{$|A|=|B|=\aleph_0$ and for every $n$ there is $k$ so
that for all $a\in A, b\in B$, if $a,b>k$ then $a/b>n$ or $b/a>n$.}

\begin{lemma} \label{noem} Suppose $A,B\su \N$  satisfy $(*)$.
Then for every infinite set $X\su D_A$ and every function $f:X\to B$,
for every $n>1$ there are distinct $x,y\in X$ so that either
$d(f(x),f(y))/d(x,y)>n$ or $d(f(x),f(y))/d(x,y)<1/n$.
\end{lemma}

\begin{proof}Let $x,y\in X$ be chosen so that $d(x,y)>0$ is sufficiently small.
\end{proof}

The lemma assures a strong form of bi-Lipschitz incomparability:
no infinite subset of one of the spaces $D_A, D_B$ can be
bi-Lipschitz embedable into the other space, if $A,B$ satisfy
$(*)$.

Let $\mathcal F$ be an almost disjoint family over $\N$, namely, a
family of infinite subsets of $\N$ with finite pairwise
intersections. Replacing each $A\in \mathcal F$ by $\{n^2:n\in
A\}$, one obtains a family of sets that pairwise satisfy
condition $(*)$. Since there is an almost disjoint family of size
$2^{\aleph_0}$ over $\N$, there is a family $\mathcal F$ of
subsets of $\N$ whose members satisfy $(*)$ pairwise, and
therefore $\{D_A:A\in \mathcal F\}$ is a family of pairwise
bi-Lipschitz incomparable subspaces of $(2^\N,d_3)$ of size
$2^{\aleph_0}$.

\begin{theorem}
Suppose $X$ is a separable metric space and that for some $K\ge
1$, for every open ball $u$ in $X$ there is a (nonempty) open
subset of $(2^\om,d_3)$ which is $K$-bi-Lipschitz embedable into
$u$. Then there are $2^{\aleph_0}$ pairwise bi-Lipschitz
incomparable $\aleph_1$-dense subsets of $X$.

In particular, in every separable Hilbert space (finite or infinite
dimensional) and in $\U$ there are $2^{\aleph_0}$ pairwise
bi-Lipschitz incomparable $\aleph_1$-dense subsets.
\end{theorem}

\begin{proof}
Fix a family $\{D_\a:\a < 2^{\aleph_0}\}$ of pairwise bi-Lipschitz
incomparable perfect subspaces of $2^\N$. For each $\a$, fix an
$\aleph_1$-dense $D^*_\a\su D_\a$.

Let $\lng u_n:n\in \N\rng$ enumerate a basis for the topology of $X$
(say all balls of rational radi with centers in some fixed countable
dense set). For each $\a<2^{\aleph_0}$, for each $n$, fix a
$K$-bi-Lipschitz embedding of a nonempty open subset of $D^*_\a$ into
$u_n$, and call the image of the embedding $E^*_{\a,n}$.  Let
$Y_\a=\bigcup_n E^*{\a,n}$. $Y_\a$ is thus an $\aleph_1$-dense subset
of $X$.

Suppose $f:Y_\a\to Y_\b$ is any function, $\a,\b<2^{\aleph_0}$
distinct, and $K>1$ is arbitrary. There is $l\in \N$ so that
$f^{-1}[E^*_{\b,l}]\cap E^*_{\a,0}$ is uncountable, hence
infinite. Therefore, by Lemma \ref{noem}, there are $x,y\in Y_\a$ so
that $d(f(x),f(y))/d(x,y)$ violates $n$-bi-Lipschitz. We conclude that
$Y_\a$, $Y_\b$ are bi-Lipschitz incomparable.

The second part of the theorem follows now from the fact that
every separable metric space embeds isometrically into $\U$ and
from the observation above that $(2^\N,d_3)$ has a bi-Lipschitz
embedding into $\R$.
\end{proof}

\section {Consistency results for separable spaces on $\aleph_1<2^{\aleph_0}$}
\label{forcing}

Let $(\mathcal M^{sep}_{\aleph_1}, \le)$ denote the set  of all
(isometry types of) separable metric spaces whose cardinality is
$\aleph_1$, quasi-ordered by almost isometric embeddability and let
$(\mathcal M_{\aleph_1},\le)$ denote the set of all (isometry types
of) metric spaces whose cardinality is $\aleph_1$, quasi-ordered
similarly.

Let $\cf(\mathcal M_{\aleph_1},\le)$ denote the \emph{cofinality} of
this quasi-ordered set: the least cardinality of $D\su \mathcal
M^{sep}_{\aleph_1}$ with the property that for every $M\in \mathcal
M^{sep}_{\aleph_1}$ there is $N\in D$ so that $M\le N$. The statement ``$\cf
(\mathcal M^{sep}_{\aleph_1},\le)=1$'' means that there is a single
$\aleph_1$-dense subset of $\U$ in which every $\aleph_1$-dense subset
of $\U$ is almost isometrically embedded, or, equivalently, that there
is a \emph{universal} separable metric space of size $\aleph_1$ for
almost isometric embeddings.

In the previous section it was shown that there are $2^{\aleph_0}$
pairwise incomparable elements --- an anti-chain --- in this quasi
ordering, with each of  elements being an $\aleph_1$-dense
subset of $\U$. This in itself does not rule-out the possibility that
a universal separable metric space of size $\aleph_1$ exists for
almost isometric embeddings.  In fact, if CH holds, $\U$ itself is
such a set.

 What can one expect if CH fails?  There has been a fairly extensive
study of the problem of universality in $\aleph_1<2^{\aleph_0}$. On
the one hand, it is fairly routine to produce models in which there
are no universal objects in cardinality $\aleph_1$ in every reasonable
class of structures (linear orders, graphs, etc) \cite{mekler,
KjSh:409} and here too it is easy to find models in which neither a
separable metric space nor a general metric space exist in cardinality
$\aleph_1$ (see below).

On the other hand, it has been shown that universal linear
orderings may exists in $\aleph_1<2^{\aleph_0}$ \cite{Sh:100},
that universal graphs may exist in a prescribed regular
$\l<2^{\aleph_0}$ \cite{Sh:gr1, Sh:gr2} and more generally, that
uninversal relational theories with certain amalgamation
properties may have universal models in regular uncountable
$\l<2^{\aleph_0}$ \cite{mekler}.

We do not know whether it is consistent to have a universal separable
metric space of size $\aleph_1<2^{\aleph_0}$ for almost isometric
embeddings, but we shall prove a weaker statement: the consistency
that $\aleph_1< 2^{\aleph_0}$ and that $\cf (\mathcal
M^{sep}_{\aleph_1},\le)$ is smaller than the continuum.

We begin by relating separable to nonseparable spaces:

\begin{theorem}
 $\cf (\mathcal
M^{sep}_{\le \aleph_n},\le)\le \cf(\mathcal M_{\aleph_n},\le)$ for all
$n$.
\end{theorem}

\begin{proof}Suppose $M=(\om_1,d)$ is a metric space. For every
  ordinal $\a<\om_1$ denote the closure of $\a$ in $(M,d)$ by $X_\a$
  is a separable space and is therefore isometric to some $Y_\a\su
  \U$. Let $N(M)=\bigcup_{\a<\om_1} Y_\a$. $N(M)$ is a subspace of $\U$
  whose cardinality is $\le \aleph_1$.

Suppose $X$ is any separable subspace of $M$, and fix a countable
dense $A\su X$. There is some $\a<\om_1$ so that $A\su \a$, hence
$X\su X_\a$ and is thus isometrically embedded in $Y_\a\su N$. In
other words, there is a single subspace of $\U$ into which all
\emph{separable} subspaces of $M$ are isometrically embedded.

Suppose now that $\cf(\mathcal M_{\aleph_1},\le)=\k$ and fix $D\su
\mathcal M_{\aleph_1}$ of cardinality $\k$ so that for all $N\in
\mathcal M_{\aleph_1}$ there is $M\in D$ so that $N\le D$. For each
$M\in D$ let $N(M)\su \U$ be chosen as above. We claim that
$\{ N(M): M\in D\}$ demonstrates that $\cf
(\mathcal M^{sep}_{\le\aleph_1},\le)\le \k$. Suppose that $X$ is a separable
metric space of cardinality $\aleph_1$. Then $X$ is almost
isometrically embedded into $M$ for some $M\in D$. Since every
separable subspace of $M$ is isometric to a subspace of $N(M)$, it
follows that $X$ is almost isometrically embedded in $N(M)$.

Simple induction on $n$ shows that for every $n$ there is a collection
$\mathcal F_n$ of $\aleph_n$ many countable subsets of $\om_n$ with
the property that every countable subset of $\om_n$ is contained in
one of them. Working with $\mathcal F_n$ instead of the collection of
initial segments of $\om_1$ gives that
$\cf(\mathcal M^{sep}_{\le\aleph_n},\le)\le \cf(\mathcal
M_{\aleph_n},\le)$.
\end{proof}

\begin{theorem}
After adding $\l\ge\aleph_2$ Cohen reals to a universe $V$ of set theory,
$\cf(\mathcal M^{sep}_{\aleph_1},\le)\ge \l$.
\end{theorem}

\begin{proof}
View adding $\l$ Cohen reals as an iteration. Let $\theta<\l$ be
given. For any family $\{A_\a:\a<\theta\}$ of $\aleph_1$-dense
subsets of $\U$ in the extension it may be assumed, by using
$\theta$ of the Cohen reals, that $A_\a\in V$ for all $\a<\theta$.
Let  $X=\Q\cup \{r_i:i<\om_1\}$ be a metric subspace of $\R$,
where each $r_i$, for $i<\om_1$ is one of the Cohen reals. We
argue that $X$ cannot be almost isometrically embedded into any
$A_\a$. Suppose to the contrary that $f:X\to A_\a$ is a
bi-Lipschitz embedding. By using countably many of the Cohen
reals, we may assume that $f\rest Q\in V$. If $f(r_0)\in A_\a$
and $f\rest Q$ are both in $V$, so is $r_0$ --- contradiction.
\end{proof}

For the next consistency result we need the following consistency
result:

\begin{theorem}[\cite{baumad}] For every regular $\l>\aleph_0$ and
  regular $\theta>\l^+$ there is a model $V$ of set theory in which
  $2^{\aleph_0} \ge \theta$ and there is a family $\{A_\a:\a<\theta\}$
  of subsets of $\l$, each $A_\a$ of cardinality $\l$ and $|A_\a\cap
  A_\b|<\aleph_0$ for all $\a<\b<\theta$.
\end{theorem}

We now state and prove the consistency for $\l=\aleph_1$, for
simplicity. Then we extend it to a general regular $\l>\aleph_0$.

\begin{theorem}
It is consistent that $2^{\aleph_0}=\aleph_3$ and that there are
$\aleph_2$ separable metric spaces on $\om_1$ such that every
separable metric space is almost isometrically embedded into one of
them.
\end{theorem}

The model of set theory which demonstrates this consistency is
obtained as a forcing extension of a ground model which satisfies
$2^{\aleph_0}=\aleph_3$ and there are $\aleph_3$
$\aleph_1$-subsets of $\aleph_1$ with finite pairwise
intersections. Such a model exists by \cite{baumad}. Then the
forcing extension is obtained via a ccc finite support iteration
of length $\aleph_2$. In each step $\z<\om_2$ a single new
separable metric space $M_\z$ of cardinality $\aleph_1$ is forced
together with almost isometric embeddings of all $\aleph_1$-dense
subsets of $\U$ that $V_\a$ knows.  At the end of the iteration,
every $\aleph_1$ dense subset is almost isometrically embedded
into one of the spaces $M_\z$ that were forced.


Let $\{A_\a:\a<\om_3\}$ be a collection of subsets of $\om_1$, each or
cardinality $\aleph_1$ and for all $\a<\b<\om_3$ it holds that
$A_\a\cap A_\b$ is finite. For each $\a<\om_3$ fix a partition
$A_\a=\bigcup_{i<\om_1} A_{\a,i}$ to $\aleph_1$ parts, each of
cardinality $\aleph_0$.

Fix an enumeration $\lng d_\a:\a<\om_3\rng$ of all metrics $d$ on
$\om_1$ with respect to which $(\om_1,d)$ is a separable metric
space and every interval $(\a,\a+\om)$ is dense in it. Since
every metric space of cardinality $\om_1$ can be well ordered in
ordertype $\om_1$ so that every interval $(\a,\a+\om)$ is a dense
set, this list contains $\om_3$ isometric copies of every
separable metric space of cardinality $\aleph_1$.

We define now the forcing notion $Q$.
A condition $p\in q$ is an ordered quintuple $p=\lng w^p,u^p,d^p,\ov
f^p, \ov \e^p\rng$ where:
\begin{enumerate}
\item $w^p$ is a finite subset of $\om_3$ (intuitively --- the set of
  metric spaces $(\om_1,d_\a)$ which the condition handles)

\item $u^p\su
\om_2$ is finite and $d^p$ is a metric over $u^p$. $(u^p,d^p)$ is a
finite approximation to the space $M=(\om_1,d)$ which $Q$ introduces.

\item $\ov\e^p=\lng \e^p_\a:\a\in w^p\rng$ is a sequence of rational
  numbers from $(0,1)$.

\item $\ov f=\lng f_\a:\a\in w^p\rng$ is a sequence of finite function
 $f_\a:(\om_1,d_\a)\to (u^p,d^p)$ that satisfy:
\begin{enumerate}
\item $f^p_\a(i)\in
 A_{\a,i}$ for each $i\in \dom f^p_\a$;
\item each $f^p_\a$ is
 $(1+\e^p)$-bi-Lipschitz.
\end{enumerate}
\end{enumerate}

The order relation is: $p\le q$ ($q$ extends $p$) iff $w^p\su w^q$,
$(u^p,d^p)$ is a subspace of $(u^q,d^q)$, and for all $\a\in
w^p$,  $\e^p=\e^q$ and $f_\a^p\su f_\a^q$.

Informally, a condition $p$ provides finite approximations of
$(1+\e^p_\a)$-bi-Lipschitz embeddings of $(\om_1,d_\a)$ for finitely
many $\a<\om_3$ into a finite space $(u^p,d^p)$ which approximates
$(\om_1,d)$.

\begin{lemma}[Density]
For every $p\in Q$:
\begin{enumerate}
\item For every $j\in \om_1\sm u^p$ and a metric type $t$ over
  $(u^p,d^p)$ there is a condition $q\ge p$ so
  that $j\in u^q$ and $j$ realizes $t$ over $u^p$ in $(u^q,d^q)$.
\item For every  $\a\in \om_3\sm w^p$ and $\d>0$ there is a
 condition $q\ge p$ so that $\a\in w^q$ and $\e^p_\a<\d$.
\item For every  $\a\in w^p$ and $i\in \om_1\sm \dom f^p_\a$ there is
 a condition $q\ge p$ so that $i\in \dom f^q_\a$
\end{enumerate}
\end{lemma}

\begin{proof}
To prove (1) simply extend $(u^p,d^p)$ to a metric space
$(u^q,d^q)$ which contains $j$ and in which $j$ realizes $t$ over
$u^p$, leaving everything else in $p$ unchanged.

For (2) define $w^q=w^p\cup \{\a\}$ and
$\e^p_\a=\e$ for a rational $\e<\d$.

For (3) suppose $i\notin \dom f^p_\a$. Fix some $x\in A_{\a,i}\sm
u^p$ and fix some $1<\l< 1+\e^p_\a$ so that $d_\a(j,k)/\l\le
d^t(j,k)\le \l d_\a(j,k)$ for all distinct $j,j\in \dom f^p_\a$.
Let $d^*(x,f^p_\a(j))=r_j $ for some rational $\l d_\a(i,j)\le
r_j < (1+\e^p_\a) d_\a(i,j)$, and let $d^*(x,x)=0$. Let $d$ be
the shortest path metric obtained from $d^t\cup d^*$. This is
obviously a rational metric and as in the proof of Fact \ref{sp},
it follows that this metric extends $d^t$ and that $f^p_\a\cup
\{(i,x)\}$ is $(1+\e^p_t)$-bi-Lipschitz into $u^t\cup \{x\}$ with
the extended metric.
\end{proof}

\begin{lemma}
In $V^Q$ there is  a separable metric space $M=(\om_1,d)$ so that for every
 separable metric space $(\om_1,d')\in V$  and $\d>0$ there is a
 $(1+\d)$-bi-Lipschitz embedding in $V^Q$ of $(\om_1,d')$ into $M$.
\end{lemma}

\begin{proof}
Let $d=\bigcup_{p\in G}d^p$ where $G$ is a $V$-generic filter of $Q$. By
(1) in the density Lemma, $d$ is a metric on $\om_1$ and, furthermore, for
every given $i<\om_1$ the interval $(i,i+\om)$ is dense in
$(\om_1,d)$.

Given any metric space $(\om_1,d')$ and a condition $p$, find some
$\a\in\om_3\sm w^p$ so that $(\om_1,d_\a)$ is isometric to
$(\om_1,d')$ and apply (2) to find a stronger condition $q$ so
that $\a\in w^q$ and $\e^q<\d$. For every $i<\om_1$ there is a
stronger condition $q'$ so that $i\in \dom f^{q'}_\a$. Thus, the
set of conditions which force a partial $(1+\d)$-bi-Lipschitz
embedding from $(\om_1,\d_\a)$ into $(\om_1,d)$ which includes a
prescribed $i<\om_1$ in its domain is dense; therefore $Q$ forces
a $(1+\d)$-bi-Lipschitz embedding of $(\om_1,d')$ into
$(\om_1,d)$.
\end{proof}

\begin{lemma}
Every antichain in $Q$ is countable.
\end{lemma}

\begin{proof}
Before plunging into the details, let us sketch shortly the main
idea of the proof.  While combining two arbitrary finite
$\l$-bi-Lipschitz embeddings into a single one is not generally
possible, when the domains are sufficiently ``close'' to each
other, that is, have very small Hausdorff distance, it is
possible. The main point in the proof is that separability of the
spaces $(\om_1,\d_\a)$ for $\a<\om_3$ implies that among any
$\aleph_1$ finite disjoint subsets of $(\om_1,d_\a)$ there are
two which are very small perturbations of each other, and
therefore have a common extension.

Suppose $\{p_\z:\a<\om_1\}\su Q$ is a set of
conditions. Applying the $\Delta$-system lemma  the
pigeon hole principle a few times, we may assume (after replacing the
set of conditions by a subset and re-enumerating) that:

\begin{enumerate}
\item $|u^{p_\z}|=n$ for all $\z<\om_1$ for some fixed $n$ and
  $\{u_\z:\z<\om_1\}$ is a $\Delta$-system with root $u$. Denote
  $u^{p_\z}=u\cup u_\z$ (so $\z<\xi<\om_1\imply u_\z\cap u_\xi=\emptyset$).
\item $d^{p_\z}\rest u=d$ for some fixed (rational) metric
  $d$.
\item Denote by $g_{\z,\xi}$  the order preserving map from
  $u_\z$ onto $u_\xi$. Then  $\id_u\cup g_{\z,\xi}$ is an
  isometry between $u^{p_\z}$ and $u^{p_\xi}$. This means
  that $u_\z$ and $u_\xi$ are isometric and  that for every $x\in
  u_\z$ and $y\in u$, $d^{p_\z}(x,y)=d^{p_\xi}(g_{\z,\xi}(x),y)$.
\item $\{w^{p_\z}:\z<\om_1\}$ is a $\Delta$-system with root $w$.
\item $\e^{p_\z}_\a=\e_\a$ for some fixed $\e_\a$ for every $\a\in w$ and
  $\z<\om_1$.
\item For every $\a\in w$, $\{\dom f^{p_\z}_\a:\z<\om_1\}$ is a
  $\Delta$ system with root $r_\a$ and $|\dom f^{p_\z}_\a|$ is
  fixed.
\item $f^{p_\z}_\a\rest r_\a$ is fixed (may be assumed since
  $f^{p_\z}_\a(i)\in A_{\a,i}$ for all $i\in r$ and $A_{\a,i}$ is
  countable).
\item  Denote $\dom f^{p_\z}_\a=r_\a\cup
  s^\z_\a$ for $\a\in w$; then  $f^{p_\z}_\a(s^\z_\a)\cap
  u=\emptyset$.
\item For all $\a\in w$ and $\z,\xi<\om_1$,
  $g_{\z,\xi}[f^{p_\z}_\a[s^\z_\a]]=f^{p_\xi}_\a[s^\xi_\a]$. Denote
  $h^{\z,\xi}_\a=f^{p_\z}_\a\circ
  h_{\z,\xi}\circ (f^{p_\xi}_\a)^{-1}$. Thus $h^{\z,\xi}_\a:s^\z_\a\to
  s^\xi_\a$ and $f^{p_\xi}_\a(h^{\z,\xi}_\a(x))=f^{p_\z}_\a(i)$ for
  all $i\in s^z_\a$.
\item For all $\a,\b\in w$ and $\z<\om_1$ it holds that $\ran(
  f^\z_\a)\cap \ran (f^\z_\b)\su u$; here we use the fact that $\ran
  f^\z_\a\su A_\a, \ran (f^\z_\b) \su A_\b$ and $|A_\a\cap
  A_\b|<\aleph_0$.
\end{enumerate}

Every $(\om_1,d_\a)$ is separable, and in a separable
metric space every subset of size $\aleph_1$ may be thinned out to a
subset of the same size in which each points is a point of complete
accumulation of the set, that is, every neighborhood of a point from
the set contains $\aleph_1$ points from the set. Therefore we may also
assume:

\begin{enumerate}
\item[(11)] for every $\a\in w$ and $\z<\om$, every $i\in s^\z_\a$ is a
point of complete accumulation in $(\om_1,d_\a)$ of
$\{h^{\z,\xi}_\a(i):\xi<\om_1\}$.
\end{enumerate}

We shall find now two conditions $p_\z,p_\xi$, $\z<\xi<\om_1$ and a
condition $t\in Q$ which extends both $p_\z$ and $p_\xi$.

Fix some $\z<\om_1$ and define $\d_0=\min\{d^{p_\z}(x,y):x\not=y\in
u\cup u_\z\}$.
Next, for $\a\in w$ let $i,j\in \dom
f^{p_\z}_\a$ be any pair of points. It holds that
\begin{equation}\label{eq1}
d_\a(i,j)/(1+\e^{p_\z}_\z)<
d^{p_\z}(f^{p_\z}_\z(i),f^{p_\z}_\z(j))<(1+\e^{p_\z}_\z)d(i,j)
\end{equation}

Which, denoting  for simplicity $\l=(1+\e^{p_\z}_\a)$,
$a=d_\a(i,j)$ and $b=d^{p_\z}(f^{p_\z}_\a(i),f^{p_\z}_\z(j))$,
is  re-written as
\begin{equation}\label{eq2}
a/\l<b<\l a
\end{equation}

There is a sufficiently small $\d>0$, depending on $a$ and $b$, such
that whenever $|b-b'|<\d$ and $|a-a'|<\d$, it holds that

\begin{equation}\label{eq3}
a'/\l<b'<\l a'
\end{equation}



Since $w$ and each $\dom f^{p_\z}_\a$ for $\a\in w$ are finite,
we can fix $\d_1>0$ which is sufficiently small for (\ref{eq3})
to hold for all $\a\in w$ and $i,j\in \dom f^{p_\z}_\a$.

 Let $\d=\min\{\d_0/100,\d_1/2 \}$.

By condition (11) above,  find $\z<\xi<\om_1$ so that for
each $\a\in w$ and $i\in s^\z_\a$ it holds that
$d_\a(i,h^{\z,\xi}_\a(i))<\d$.

Let $u^t=u\cup u_\z\cup u_\xi$. We define now a metric $d^t$ on $u^t$
as follows. First, for each $\a\in w$ and $x=f^{p_\z}_\a(i)\in
f^{p_\z}_\a[s^\z_\a]$ let $y=g_{\z,\xi}(x)$ and define
$d^*(x,y)=r$, $r$  a rational number which
satisfies $1/(1+\e^p_\a)d_\a(i,j) < r < (1+\e^p_\a)d_\a(i,j)$, where
$j=h^{\z,\xi}_\a(i)$. Since for $\a,\b\in w$ the sets
$f^{p_\z}_\a[s^\z_\a]$ and $ f^{p_\z}_\a[s^\xi_\b]$ are disjoint, and
similarly for $\xi$, $d^*$ is well defined, namely, at most one $\a$
is involved in defining the distance $r$.  In fact, any two $d^*$
edges are vertex disjoint.

Now  $(u\cup u_\z\cup u_\xi,d\cup d^*)$  is a connected weighted
graph.  Let $d^t$ be the shortest-path metric on $u\cup u_\z\cup u_\xi$
obtained from $d^{p_\z}\cup d^{p_\xi}\cup d^*$. It is obviously a
rational metric, as all distances in $d^{p_\z}\cup d^{p_\xi}\cup d^*$
are rational.

Let us verify that $d^t$ extends $d^{p_\z}\cup d^{p_\xi}\cup
d^*$. Suppose $x\in u_\z,y\in u_\xi$ and $d^*(x,y)$ is defined. Any
path from $x$ to $y$ other than $(x,y)$ must contain some edge with a
distance in $(u^{p_\z},d^{p_\z})$, and all those distances are much
larger than $d^*(x,y)$, so $(x,y)$ is the shortest path from $x$ to
$u$. Suppose now that $x,y\in u^{p_\z}$. So $(x,y)$ is the shortest
path among all paths that lie in $u^{p_\z}$ and the path of minimal
length from $x$ to $y$ among all paths that contain at least one $d^*$
edge is necessarily the path $(x,x',y',y)$ where $x'=g_{\z,\xi}(x)$
and $y'=g_{\z,\xi}(y)$ (since $(u^{p_\xi},d^{p_\xi})$ is a metric
space).  Since $d^{p_\z}(x,y)=d^{p_\xi}(x',y')$ by condition (3), the
length of this path is larger than $d^{p_\z}(x,y)$.

Let $f^t_\a=f^{p_\z}_\z\cup f^{p_\xi}_\z$, where we formally take
$f^p_\a$ to be the empty function if $\a\notin w^p$.

Let $w^t=w^{p_\z}\cup w^{p_\xi}$ and let
$\e^t_\a=\max\{\e^{p_\z}_\a,\e^{p_\xi}_\a\}$ where $\e^p_\a$ is taken
as $0$ if $\a\notin w^p$ (recall that $\e^{p_\z}_\a=\e^{p_\xi}_\a$ for
$\a\in w$).

Now $t$ is defined, and extends both $p_\z$ and $p_\xi$; one only
needs to verify that $t\in Q$. For that, we need to verify that each
$f^t_\a$ is $(1+\e^t_\a)$-bi-Lipschitz. For $\a\notin w$ this is
trivial. Suppose $\a\in w$ and $i,h\in \dom f^t_\a=r_\a\cup s_\z\cup
s_\xi$. The only case to check is when $i\in s^\z_\a$ and $j\in
s^\xi_\a$. If  $h^{\z,\xi}_\a(i)=j$ then this is
taken care of  by the choice of
$d^t(f^{p_\z}_\a(i),f^{p_\xi}_\a(j))=d^*(f^{p_\z}_\a(i),f^{p_\xi}_\a(j))$.

We are left with the  main case:  $j\not=h^{\z,\xi}_\a(i)$.  Let
$x=f^{p_\z}_\a(i)$, $y=f^{p_\xi}_\a(j)$ and denote $a'=d_\a(i,j)$.
Let $j'\in x^\z_\a$ be such  that $h^{\z,\xi}_\a(j')=j$ and let
$y'=g_{\z,\xi}^{-1}(y)$ (so $f^{p_\z}_\z(j')=y'$).

Denote $b=d^t(x,y')$, $a=d_\a(i,j')$.  We have that $a/\l<b<\l a$, and
need to prove  $a'/\l <b'<\l b'$.

It holds that $d_\a(j,j')<\d\le\d_1/2$, hence $d^t(y,y')<\d_1$.  By
the triangle inequality in $u^t$, we have that $|b-b'|<\d_1$. On the
other hand, by the triangle inequality in $(\om_1,d_\a)$, we have that
$|a'-a|<\d<\d_1$. Thus, by the choice of $\d_1$ so that (\ref{eq3})
holds if (\ref{eq2}) holds and $|a-a'|<\d_1$, $|b-b'|<\d_1$, we have
that $a'/\l<b'<\l a'$, as required.
\end{proof}

 Let $P=\lng P_\b,Q_\b:\b<\om_2\rng$ be a finite support iteration of
length $\om_2$ in which each factor $Q_\b$ is the forcing notion we
defined above, in $V^{P_\b}$. Since each $Q_\b$ satisfies the ccc, the
whole iteration satisfies the ccc and no cardinals or cofinalities are
collapsed in the way --- in particular the collection
$\{A_\a:\a<\om_3\}$ required for the definition of $Q$ is preserved.

Since every metric $d$ on $\om_1$ appears in some intermediate stage,
the universe $V^P$ satisfies that there is a collection of $\om_2$
separable metrics on $\om_1$ so that every separable metric on $\om_1$
is almost isometrically embedded into one of them.

There is no particular property of $\om_1$ that was required in the
proof. Also, Baumgartner's result holds for other cardinals. We have
proved then:

\begin{theorem}
Let $\l>\aleph_0$ be a regular cardinal. It is consistent that
$2^{\aleph_0}>\l^+$ and that there are $\l^+$
separable metrics on $\l$ such that every separable metric on $\l$ is
almost-isometrically embedded into one of them.
\end{theorem}

In the next section we shall see that separability is essential
for this consistency result.

\section{Non-separable spaces below the continuum}
\label{club}

Now we show that the consistency proved in the previous Section
for separable metric spaces of regular cardinality
$\l<2^{\aleph_0}$ is not possible for metric spaces in general if
$\l>\aleph_1$. Even a weaker fact fails: there cannot be fewer
than continuum metric spaces on $\l$ so that every metric space
is bi-Lipschit embedable into one of them if
$\aleph_1<\l<2^{\aleph_0}$ and $\l$ is regular.

\begin{theorem}If $\aleph_1<\l<2^{\aleph_0}$ then for every
$\k<2^{\aleph_0}$ and metric spaces $\{ (\l,d_i):i<\k\}$ there
exists an ultra-metric space $(\l,d)$ that is not bi-Lipschitz
embedable into $(\l,d_i)$ for all $i<\k$. In particular there is
no single metric space $(\om_2,d) $ into which every ultra-metric
space of cardinality $\l$ is bi-Lipschitz embedded.
\end{theorem}

\begin{proof}

Let $\l>\aleph_1$ be a regular cardinal.  Let
$S^\l_0=\{\d<\l:\cf\d=\l\}$, the stationary subset of $\l$ of
countably cofinal elements. For a regular $\l>\aleph_1$ we may fix a
club guessing sequence $\ov C=\lng c_\a:\d\in S^\l_0\rng$ \cite{CA,KjSh:409}:

\begin{enumerate}
\item $c_\d\su\d=\sup c_\d$ and $\otp c_\d=\om$ for all $\d\in S^\l_0$;
\item For every club $E\su \l$ the set $S(E)=\{\d\in S^\l_0:c_\d\su
E\}$ is stationary.
\end{enumerate}

 For each $\d\in S^\l_0$ let $\lng \a^\d_n:n<\om\rng$ be the increasing
 enumeration of $c_\d$.

 Let $(\l,d)$ be a given metric space. Let $X_\a$ denote the
 subspace $\{\b:\b<\a\}$. Let $d(\b,X_\a)$ denote the distance of $\b$
 from $x_\a$, that is, the infimum of $d(\b,\g)$ for all $\g<\a$. Thus
 $\lng d(\b,X_{\a_n}):n<\om\rng$ is a weakly decreasing sequence of
 positive real numbers.

\begin{lemma}
Suppose $f:\l\to \l$ is a bi-Lipschitz embedding of
$(\l,d_1)$ in $(\l,d_2)$ with constant $K\ge 1$. Then there is a
club $E\su \l$ such that: for all $\a\in E$ and $\b>\a$ it holds
that $f(\b)>\a$ and
\begin{equation}\label{ineq}
(2K)^{-1}d_1(\b,X_\a)\le d_2(f(\b),X_\a)\le Kd_1(\b,X_\a)
\end{equation}

\end{lemma}

\begin{proof}
Consider the structure $\mathcal M= (\l; d_1, d_2, f, \lng
P_q:q\in \Q^+\rng)$ where $P_q$ is a binary predicate so that
$\mathcal M\sat P_q(\b_1,\b_2)$ iff $d_2(\b_1,\b_2)< q$.

Let $E=\{\a<\l: M\rest X_\a\prec M\}$. Then $E\su \l$ is a
club. Suppose that $\a\in E$. Then $X_\a$ is closed under $f$ and
$f^{-1}$ and therefore $f(\b)>\a$.  The inequality
$d_2(f(\b),X_\a)\le Kd_1(\b,X_\a)$ is clear because $X_\a$ is
preserved under $f$. For the other inequality suppose that $\g\in
X_\a$ is arbitrary and let $\e:=d_2(f(\b),\g)$. Let $q>\e$ be an
arbitrary rational number. Now $\mathcal M\sat P_q(f(\b),\g)$ and
therefore, by $\mathcal M\rest X_\a\prec \mathcal M$, there exists
some $\b'\in X_\a$ so that $\mathcal M\sat P_q((f(\b'),\g)$.  Thus

\[d_2(f(\b),f(\b'))\le d_2(f(\b),\g)+d_2(\g,f(\b'))< 2q
\]

Since $K^{-1}d_1(\b,\b') \le d_2(f(\b),f(\b'))$ it follows that
$K^{-1}d_1(\b,X_\a)< 2q$, and since $q>\e$ is an arbitrary rational,
it follows that $(2 K)^{-1}d_1(\b,X_\a)\le d_2(f(\b),X_\a)$.
\end{proof}

Suppose $(\l,d_i)$ are metric spaces for $i=0,1$, that $K\ge 1$ and
that $f:\l\to \l$ is a $K$-bi-Lipschitz embedding of
$(\l,d_1)$ in $(\l,d_2)$. Let $E\su \l$ be a club as
guaranteed by the previous lemma. Suppose $\d\in S(E)$, $\b>\d$ and
$\lng \a_n:n<\om\rng$ is the increasing enumeration of $c_\d$. Let
$\e_n=d_1(\b,X_{\a_n})$ and let $\e'_n=d_2(f(\b),X_{\a_n})$. From
(\ref{ineq}) it follows that for each $n$,

\begin{equation}\label{ineq2}
 (2K^2)^{-1}\e_n/ \e_{n+1} \le \e'_n/\e'_{n+1} \le 2K^2\e_n/\e_{n+1}
\end{equation}

 For a subset $A\su
\om$, $\d\in S^2_0$, $\b>\d$ a metric $d$ over $\om_2$   and an integer
$K\ge 1$ we write $\Theta_d(\b,\d,A,K)$ if $\e_n/\e_{n+1} =1 $ or
$\e_n/\e_{n+1}>4K^4$ for all $n$ and $A=\{n:\e_n/\e_{n+1}>4K^4\}$.

Write $\Phi_d(\b,\d,A,K)$ if $A=\{n:\e_m/\e_{n+1}>2K^2\}$.  From
(\ref{ineq2}) it follows:

\begin{lemma}\label{pres}
Suppose $f:\om_2\to \om_2$ is a $K$-bi-Lipschitz embedding of
$(\om_2,d_1)$ in $(\om_2,d_2)$. Then there is a club $E\su \om_2$ such
that for all $\d\in S(E)$, $A\su \om$ and $\b>\d$ it holds that
$f(\b)>\d$ and $\Phi_{d_1}(\b,\d,A,K)\imply
\Theta_{d_2}(f(\b),\d,A,K)$.
\end{lemma}

\begin{lemma} \label{construction} For every infinite $A\su \om$  there is
an ultra-metric space $(\om_2,d)$ and a club $E\su \om_2$ so that for
all $\d\in S(E)$ and integer $K\ge 1$ there exists $\b>\d$ with
$\Phi_d(\b,\d,A,K)$.
\end{lemma}

\begin{proof} Suppose $A\su \om$ is infinite and let
$\lng a_n:n<\om\rng$ be the increasing enumeration of $A$.  Let
$(\om_2)^\om$ be the tree of all $\om$-sequences over $\om_2$ and let
$d$ be the metric so that for distinct $\eta_1,\eta_2\in (\om_2)^\om$,
$d(\eta_1,\eta_2)=1/(n+1)$, where $n$ is the length of the largest
common initial segment of $\eta_1,\eta_2$ (or, equivalently,
$d(\eta_1,\eta_2)=1/(|\eta_1\cap \eta_2|+1)$). Every
\emph{finite} sequence $t:n\to \om_2$ determines a basic clopen
ball of radius $1/(n+1)$ in $((\om_2)^\om, d)$, which is
$B_t=\{\eta\in (\om_2)^\om: t\su \eta\}$.

By induction on $\a<\om_2$ define an increasing and continuous chain
of subsets $X_\a\su(\om_2)^\om$ so that:
\begin{enumerate}
\item $|X_\a|=\om_1$ for all $\a<\om_2$
\item $X_\a\su X_{\a+1}$ and $X_\a=\bigcup_{\b<\a}$ if $\a<\om_2$ is limit.
\item For every $\nu\in X_\a$ and $k<\om$  there exists $\eta\in X_{\a+1}\sm
X_\a$ so that $\nu\rest k\su \eta$ but $B_{\nu\rest (n+1)}\cap
X_\a=\emptyset$.
\item If $\a=\d\in S^2_\a$, $\lng \a_n:n<\om\rng$ is the
increasing enumeration of $c_\d$, then for every integer $K\ge 1$ a
sequence $\nu_K$ defined as follows. Let $b_n=(4K^2)^{n+1}$.
Define an increasing sequence of finite sequences $t_n$ with
$|t_n|=b_n$ by induction on $n$ as follows. $t_0=\lng \rng$.
Suppose that $t_n$ is defined and $B_{t_n}\cap
X_{\a_{a_n}}\not=\emptyset$.  Choose $\eta\in X_{\a_{a_n}+1}\cap
B_{t_n}$ so that $B_{\eta\rest (b_{n+1})}\cap
X_{\a_{a_n}}=\emptyset$. Now let $t_{n+1}=\eta\rest (b_{n+1})$.

 Finally, let $\nu_K =\bigcup_n t_n$.  Put $\nu_k$ in $X_{\d+1}$ for
 each integer $K\ge 1$.
\end{enumerate}

There is no problem to define $X_\a$ for all $\a<\om_2$. Let
$X=\bigcup_{\a<\om_2}X_\a$. Fix a 1-1 onto function $F:X\to
\om_2$ and  let $d$ be the metric on $\om_2$ which makes $F$ an
isometry. Observe that for some club $E\su \om_2$ it holds that
$F[X_\a]=\a$ for all $\a\in E$. If $c_\d\su E$ let $\nu_K$ be the
sequence we put in $X_{\a+1}$ in clause (4) of the inductive
definition and let $\b_K=F(\nu_K)$. We leave it to the reader to
verify that $\Phi_d(\b_K,\d,A,K)$.
\end{proof}

The proof of the Theorem follows from both lemmas.  Suppose
$\k<2^{\aleph_0}$ and $d_i$ is a metric on $\om_2$ for each
$i<\k$. The set $\{A\su \om: (\exists i<\k)(\exists
\d\in S^2_0)(\exists K\in \om\sm \{0\})( \exists \b>\d)
[ \Theta_{d_i}(\b,\d,A,K)]
\}$ has cardinality $\le \aleph_2\cdot\k<2^{\aleph_0}$. Therefore there
is some infinite $A\su \om$ not in this set. By Lemma
\ref{construction} there is some ultra-metric $d$ on $\om_2$ and a club
$E\su\om_2$ so that for all $\d\in S(E) $ and integer $K\ge 1$ there
is $\b>\d$ with $\Phi_d(\b,\d,A,K)$. If for some $i<\k$ and integer
$K\ge 1$ there is a $K$-bi-Lipschitz embedding $\phi$ of $(\om_2,d)$
into $(\om_2,d_i)$ then from Lemma
\ref{pres} there is some club $E'\su\om_2$ so that
 for all $\d\in S(E)$, $\b>\d$ and infinite $A\su \om$, $\phi(\b)>\d$
 and $\Phi_d(\b,\d,A, K) \imply \Theta_{d_i}(\b,\d,A,K)$. Let
 $E''=E\cap E'$. Choose $\d\in S(E'')$. Now for some $\b>\d$ it holds
 $\Phi_d(\b,\d,A)$, hence $\Theta_{d_i}(\b,\d,A,K)$ contrary to the
 choice of $A$.
\end{proof}

Let us take a second look at the proof above. What is shown is
that subsets of $\N$ --- ``reals'' --- can be coded into a metric
space $X$ of regular cardinality $\l$ and retrieved from a larger
space $Y$ in which $X$ is bi-Lipschitz embedded. The retrieval is
done modulo some nonstationary subset of the cardinal $\l$.

For a metric space $X$ of cardinality $\l$, let us define $S(X)$
as the set of all $A\su \N$ so that: for every enumeration of $X$
of ordertype $\l$, and for every club $E\su \l$ there is some
$\d\in S^\l_0$ so that $c_\d\su E$ and for every $K>1$ in $\N$
there is some $\b>\d$ so that $\Phi(\b,\d,A,K)$. The set of those
$A$s is preserved under bi-Lipschitz embeddings.

\begin{theorem}For every regular $\l>\aleph_1$ there is an order
  preserving map $F:(\mathcal M_\l,\le_{BL}\to ([\mathcal P(\N)]^{\le
  \l}, \su)$, where $\le_{BL}$ is the quasi-ordering of bi-Lipschitz
  embeddability and $[\mathcal P(\N)]^{\le}$ is the family of subsets
  of $\mathcal P(\N)$ whose cardinality is at most $\l$. The range of
  $F$ restricted to ultra-metric spaces is cofinal in $(\mathcal
  P(\N),\mathord{\su})$.
\end{theorem}

The picture at large is as follows: the set of all distances in a
metric space is an invariant which is preserved under isometric
embeddings which forbids the existence of  universal metric space with
respect to isometries in cardinalities below the continuum. With
respect to almost isometric embeddings there is no such invariant for
countable spaces, or for separable spaces of regular cardinality below
the continuum. However, for general spaces of regular cardinality
$\aleph_2$ or higher such invariants are again available, as proved
above, and prohibit the existence of universal metric spaces.


\begin{thebibliography}{9}

\bibitem{baumgartner} James E. Baumgartner.
{\em All $\aleph \sb{1}$-dense sets of reals can be isomorphic.},
Fund. Math. 79 (1973), no. 2, 101--106.

\bibitem{baumad} James E. Baumgartner. {\em Almost-disjoint sets, the
dense set problem and the partition calculus}, Ann. Math. Logic 9
(1976), no. 4, 401--439.


\bibitem{DzSh} M. Dzamonja and S. Shelah. {\em Universal graphs at the
  successor of a singular cardinal}, J. Symbolic Logic 68 (2003),
  no. 2, 366--388.

\bibitem{Gromov} Misha Gromov. {\em Metric structures for
Riemannian and non-Riemannian spaces}, With appendices by M. Katz,
P. Pansu and S. Semmes. Translated from the French by Sean Michael
Bates. Progress in Mathematics, 152. Birkh\" user Boston, Inc.,
Boston, MA, 1999. xx+585 pp.


\bibitem{HZ}Michael Hrusak and Beatriz Zamora
Aviles. {\em Countable dense homogeneity of
product spaces"}, preprint.

\bibitem{CGK} John D. Clemens, Su Gao and Alexander Kechris.  {\em
Polish metric spaces: their classification and isometry
groups}, Bull. Symbolic Logic 7 (2001), no. 3, 361--375.

\bibitem{GK}Su Gao and Alexander Kechris. {\em On the classification
of Polish metric spaces up to isometry}, Mem. Amer. Math. Soc. 161
(2003), no. 766, viii+78 pp.

\bibitem{rep} M. Kojman. {\em Representing embeddability as set
inclusion}, J. London Math. Soc. (2) 58 (1998), no. 2,
257--270.





 \bibitem{groups} M.  Kojman and S. Shelah.
 {\em Universal abelian groups}, Israel J. Math. 92 (1995), no. 1-3,
 113--124.



\bibitem{KjSh:409} M.  Kojman and S. Shelah.
{\em Nonexistence of
  universal orders in many cardinals}, J. Symbolic Logic 57 (1992),
  no. 3, 875--891.




\bibitem{stuss} M.  Kojman and S. Shelah.  {\em The universality
  spectrum of stable unsuperstable theories}, Ann. Pure Appl. Logic 58
  (1992), no. 1, 57--72.



\bibitem{mekler} Alan H. Mekler. {\em Universal structures in power
  $\aleph_1$}, J. Symbolic Logic 55 (1990), no. 2, 466--477.

\bibitem{CA} S. Shelah. {\em Cardinal arithmetic}, Oxford Logic
Guides, 29. Oxford Science Publications. The Clarendon Press, Oxford
University Press, New York, 1994. xxxii+481 pp.

\bibitem{Sh:100} S. Shelah. {\em Independence results}, J. Symbolic
Logic 45 (1980), no. 3, 563--573.

\bibitem{Sh:gr1} S. Shelah. {\em
 On universal graphs without
 instances of CH}, Ann. Pure Appl. Logic 26 (1984), no. 1, 75--87.

\bibitem{Sh:gr2} S. Shelah. {\em Universal graphs without instances of
${\rm CH}$: revisited}, Israel J. Math. 70 (1990), no. 1,
69--81.


\bibitem{vershik} A. M.  Vershik.  {\em The universal Uryson space,
Gromov's metric triples, and random metrics on the series of natural
numbers}, (Russian) Uspekhi Mat. Nauk 53 (1998), no. 5(323), 57--64;
translation in Russian Math. Surveys 53 (1998), no. 5, 921--928

\bibitem{woodin} W. Hugh Woodin. {\em The continuum hypothesis. II},
Notices Amer. Math. Soc. 48 (2001), no. 7, 681--690.
\end{thebibliography}
\end{document}